\documentclass[letterpaper, 10 pt, conference]{ieeeconf}  

\IEEEoverridecommandlockouts                              

\usepackage{graphicx}      
\usepackage{cite}
\usepackage{xparse}
\NewDocumentCommand{\citep}{ o o m }{\cite{#3}}

\usepackage{amsmath,amssymb,mathtools}
\usepackage{units}
\usepackage{xspace}
\usepackage{subcaption}

\usepackage{tikz}
\usetikzlibrary{external}
\usetikzlibrary{positioning}
\tikzset{external/system call={pdflatex \tikzexternalcheckshellescape -halt-on-error
		-interaction=batchmode -jobname "\image" "\texsource"}}
\usepackage{xcolor}

\usepackage{acronym}
\acrodef{mm}[MM]{min-max}
\acrodef{mg}[MG]{microgrid}
\acrodef{pv}[PV]{photovoltaic}
\acrodef{ems}[EMS]{energy management system}
\acrodef{mpc}[MPC]{model predictive control}
\acrodef{der}[DER]{distributed energy resource}
\acrodef{res}[RES]{renewable energy sources}
\acrodef{mppt}[MPPT]{maximum power point tracking}
\acrodef{soc}[SOC]{state of charge}
\acrodef{iff}[iff]{if and only if}
\acrodef{tub}[TU Berlin]{Technische Unversität Berlin}

\newcommand{\hh}{\nobreak\hspace{0pt}-\hspace{0pt}}

\newcommand{\ie}{i.e.\@\xspace}

\newcommand{\eg}{e.g.\@\xspace}

\newcommand{\T}{^\mathsf{T}}

\newcommand{\Npr}{N_p}
\newcommand{\Nsi}{N_s}

\newcommand{\RR}{\mathbb{R}}
\newcommand{\N}{\mathbb{N}}
\newcommand{\Rp}{\mathbb{R}_{>0}}

\newcommand{\Rpz}{\mathbb{R}_{\geq0}}
\newcommand{\Rnz}{\mathbb{R}_{\leq0}}

\newcommand{\tmin}{\text{min}}
\newcommand{\tmax}{\text{max}}

\newcommand{\tst}{\text{s}}
\newcommand{\tth}{\text{t}}
\newcommand{\trs}{\text{r}}
\newcommand{\tld}{\text{d}}
\newcommand{\ok}{(k)}

\newcommand{\okk}{(k+1)}
\newcommand{\okm}{(k-1)}

\newcommand{\okNp}{(k+\Npr)}

\newcommand{\Ts}{\mathrm{T}_\text{s}}

\newcommand{\pmin}{p^\tmin}
\newcommand{\pmax}{p^\tmax}
\newcommand{\umin}{u^\tmin}
\newcommand{\umax}{u^\tmax}

\newcommand{\pth}{p_\tth}
\newcommand{\pst}{p_\tst}

\newcommand{\nt}{n_\tth}
\newcommand{\ns}{n_\tst}
\newcommand{\nr}{n_\trs}
\newcommand{\nd}{n_\tld}

\newcommand{\wld}{w_\tld}

\newcommand{\ptk}{\pth\ok}

\newcommand{\psk}{\pst\ok}

\newcommand{\prk}{p_\trs\ok}
\newcommand{\utk}{u_\tth\ok}
\newcommand{\usk}{u_\tst\ok}
\newcommand{\urk}{u_\trs\ok}

\newcommand{\wrk}{w_\trs\ok}
\newcommand{\wlk}{w_\tld\ok}

\newcommand{\psik}{p_{\tst,i}\ok}

\newcommand{\ptmin}{p_\tth^\tmin}
\newcommand{\ptmax}{p_\tth^\tmax}
\newcommand{\psmin}{p_\tst^\tmin}
\newcommand{\psmax}{p_\tst^\tmax}
\newcommand{\psminn}{\bar{p}_\tst^\tmin}
\newcommand{\psmaxx}{\bar{p}_\tst^\tmax}
\newcommand{\psmink}{\psminn\ok}
\newcommand{\psmaxk}{\psmaxx\ok}

\newcommand{\prmax}{p_\trs^\tmax}
\newcommand{\prmaxk}{\bar{p}_\trs^\tmax\ok}

\newcommand{\xsmin}{x^\tmin}
\newcommand{\xsmax}{x^\tmax}

\newcommand{\deltath}{\delta_\tth}
\newcommand{\deltati}{\delta_{\tth,i}}
\newcommand{\deltatk}{\deltath\ok}
\newcommand{\deltatkm}{\deltath\okm}

\newcommand{\deltatik}{\deltati\ok}

\newcommand{\pVec}{\mathbf{p}}
\newcommand{\uVec}{\mathbf{u}}
\newcommand{\wVec}{\mathbf{w}}

\newcommand{\wminVec}{\wVec^\tmin}
\newcommand{\wmaxVec}{\wVec^\tmax}

\newcommand{\deltathVec}{\boldsymbol{\delta}_\tth}

\newcommand{\rhok}{\rho\ok}

\newcommand{\Cth}{C_\mathrm{t}}
\newcommand{\Cst}{C_\mathrm{s}}
\newcommand{\CthOn}{C_{\text{on}}}
\newcommand{\CthSw}{C_{\text{sw}}}

\newcommand{\sat}{\operatorname{sat}}

\newtheorem{problem}{Problem}
\newtheorem{prop}{Proposition}
\providecommand{\proof}[1][]{\ifthenelse{\equal{#1}{}}{\noindent\hspace{2em}{\itshape Proof: }}{\noindent\hspace{2em}{\itshape Proof #1: }}\@\xspace}

\usepackage{booktabs}
\usepackage{tabularx}
\usepackage{pgfplots}
\pgfplotsset{compat=1.18}

\begin{document}


\title{\LARGE \bf
    Saturation-aware robust optimal operation control of microgrids based on minimum-regret optimization
}

\author{Ujjwal Pratap$^{1}$, Steffen Hofmann$^{2}$ and Christian A. Hans$^{1}$%
\thanks{$^{1}$Automation and Sensorics in Networked Systems, University of Kassel, Germany,
{\tt\small \{ujjwal.pratap, hans\}@uni-kassel.de}}
\thanks{$^{2}$Control Systems Group, TU Berlin, Germany, {\tt\small steffen.hofmann@tu-berlin.de}}%
\thanks{This work was supported by BMWE through the RESUME project (03EI6042A) and IntelliEMEs project (03EN1122A).}
}

\maketitle

\begin{abstract}
This paper studies robust optimal operation control problems for microgrids with a high share of renewable energy sources. 
The main goal is to ensure an optimal operation in the presence of a wide range of scenarios of uncertain infeed of renewable sources and uncertain load demand. 
We formally state a minimum-regret robust \ac{mpc} problem and address it by making effective use of a hierarchical microgrid control structure. 
In detail, we consider an enhanced primary control layer composed of droop control and an autonomous limitation of power and energy. 
We prove that this enables us to use constant power setpoints to achieve an optimal operation under certain conditions. 
To obtain a tractable controller, we then combine the abovementioned constant saturation-aware setpoints with an energy management system, which solves a robust unit commitment problem within a model predictive control framework. 
In a case study, we finally demonstrate the viability of the control design.
\end{abstract}

\begin{keywords}
Operation control, energy management system, hierarchical control, microgrids.
\end{keywords}


\section{Introduction}
Energy production is shifting from conventional generators to inverter\hh interfaced \ac{res} worldwide.
Many \ac{res} work in grid\hh following mode and inject the available power into the grid.
This can lead to situations where there is insufficient capacity of units that perform frequency and voltage regulation to compensate power imbalance \citep{JAF2012}.
This is especially true for \acp{mg} with high share of \ac{res} \citep{GVM2011}.
It has been proposed by \citep{KhaWuLi2024} to use \ac{res} that  are interfaced using grid\hh forming inverters or inverters with advanced control schemes.
However, unlike conventional sources, \ac{res} such as \ac{pv} plants and wind turbines are intermittent and subject to uncertainties.
In order to participate in frequency regulation and power sharing, a unit has to be able to vary its power output and thus has to be operated sufficiently far from its weather-dependent physically available power \citep{GVM2011}.
In this paper, we assume at least one grid\hh forming unit is present.
Moreover, we rely on units with frequency control capabilities that restrict their output power according to physical limitations.
Implementation of this behavior on the inverter level and related issues such as stability have been widely studied, \eg \cite{DuLas2017,YazFerDavSha2020,GaoZhoAnvBla2025} and in the review paper \cite{BaeChaLuJohSeo2024}.

To ensure feasible control actions and reduce operating costs in an \ac{mg}, the \acf{ems} computes and provides power setpoints to the units. 
\ac{ems} design for microgrids has been studied extensively in the literature. This includes hierarchical control as well as \ac{ems} design for scheduling and dispatch under uncertainty \citep{OMA2014}, both of which are also the main focus of this paper. 
In this context, \ac{mpc} is among the popular approaches. It is a receding\hh horizon optimization technique that uses forecasts and can enforce unit constraints \citep{ParRikGli2014,HNR2014}. 
Most \ac{ems} approaches for \ac{mg}s \citep{OMA2014,ParRikGli2014,SZO2016} do not explicitly consider droop control in the lower control layers or the limiting behavior of units, which can lead to a mismatch between predicted and actual behavior in presence of disturbances and uncertain conditions.

The authors in \citep{HSH2021} show how saturating droop\hh control can increase overall performance in conjunction with a min\hh max \ac{mpc}\hh based \ac{ems}.
In particular, saturating droop increases the set of feasible controls, which in turn can result in lower closed\hh loop costs. However, a major drawback of min\hh max \ac{mpc} is that only worst\hh case performance is optimized.
This can lead to wasted potential as the actual realizations of uncertain available renewable power and load demand are typically not close to the worst\hh case.
This motivates the present paper which aims to achieve near\hh optimal performance under all disturbance\hh scenarios and at the same time simplify the optimization problem.

In what follows, we extend the approach in \cite{HSH2021}. In detail, we look beyond worst-case scenarios of renewable generation and load demand, while considering droop\hh controlled units and autonomous power and energy limiting features. 
Unlike \cite{HSH2021}, which optimizes only the worst-case cost, we formulate a regret-based robust \ac{mpc} problem which considers all possible disturbance scenarios.
For this problem, we provide a simplified, not necessarily optimal, solution based on constant setpoints and a unit commitment problem.

The remainder of the paper is organized as follows. 
Section~\ref{sec:microgrid:model} recalls the saturation-aware microgrid model.
Section~\ref{sec:problemFormulation} formulates the minimum-regret robust operation problem for uncertain renewable generation and load demand. 
Section~\ref{sec:ruleBasedSolution} presents the proposed solution approach, first deriving constant saturation-based control setpoints under suitable assumptions and then combining it with a robust unit-commitment problem. 
Finally, Section~\ref{sec:caseStudy} evaluates the proposed controller in a case study and compares it with a benchmark and the controller from \cite{HSH2021}.


\section{Saturation-aware microgrid Model}
\label{sec:microgrid:model}

In this section, we describe the saturation-aware islanded \ac{mg} model based on \citep{HSH2021}. 
It includes renewable power sharing and considers saturation at the physical limits of the units.
The model incorporates three types of units, which are conventional units, \eg, thermal generators, storage units, \ac{res}, \eg, wind farms and \ac{pv} power plants, as well as loads. 
The model, which will be employed in the \ac{ems} design later on, describes the power output of each unit and accounts for power sharing and saturations at unit bounds.
In addition, uncertain renewable generation and load demand are considered since unit power outputs may differ from the power setpoints provided by the \ac{ems}.

In what follows, we assume that the primary control ensures a stable operation of the \ac{mg} and that there is no secondary control. Moreover, storage losses are assumed to be negligible. 
Finally, at the \ac{ems} time scale, switching times of the conventional units are assumed to be negligible.

\subsection{Preliminaries}
$\mathbb{R},$ $\Rnz$ and $\mathbb{R}_{\geq 0},$ denote the sets of real numbers, non-positive real numbers, and non-negative real numbers, respectively. Moreover, $\mathbb{N}$ denotes the set of positive integers.
$\min(x, y)$ and $\max(x, y)$ for the vectors $x$ and $y$ provides the element\hh wise minimum and maximum respectively. 
Given matrices $A \in \mathbb{R}^{n \times m}$ and $B \in \mathbb{R}^{n \times m}$, then $A\leq B$ is true if and only if it is true for all elements at matching positions.
$\mathbf{1}$ is a vector of all ones and $\mathbf{0}$ a vector of all zeros, both of appropriate size. \text{diag}(x) is a diagonal matrix, with the elements of the vector $x$ on its diagonal.

Consider a variable $z \in \mathbb{R},$ with $z_{\min}\leq z \leq z_{\max}$. Then a saturation of $z$ is defined as
\begin{equation}
  \label{eq:saturationOperator}
    \text{sat}(z_{\text{min}},z,z_{\text{max}}) = 
    \begin{cases} 
        z_{\text{min}}, & \text{if } z < z_{\text{min}},\\
        z, & \text{if } z_{\text{min}} \leq z \leq z_{\text{max}}, \\ 
        z_{\text{max}}, & \text{if } z > z_{\text{max}}.
    \end{cases}
\end{equation}
Similarly, for $z_{\min}, z, z_{\max} \in \mathbb{R}^{\text n},$ the saturation operator $\text{sat}(.,.,.)$ is defined element-wise.

\subsection{Model}
In islanded \ac{mg}s, generation and consumption must be balanced at every time instant $k \in \mathbb{N}$, \ie,
\begin{equation}
\label{eq:power:balance}
 \mathbf{1}\T\ptk + \mathbf{1}\T\psk + \mathbf{1}\T\prk + \mathbf{1}\T\wlk = 0,
\end{equation}
where $p_\tth\ok \in \Rpz^{\nt}$, $p_\tst\ok \in \RR^{\ns},$ $p_\trs\ok \in \Rpz^{\nr}$ denote the power outputs of conventional units, storage units and \acs{res}, respectively, and $w_\tld\ok \in \Rnz^{\nd}$ is the load with $\nt, \ns, \nr, \nd\in \mathbb{N}$ being the corresponding numbers of units.

Since load and renewable generation are uncertain, grid-forming units must adjust their power outputs so that \eqref{eq:power:balance} is satisfied. 
This can be done using droop control. 
Since the \ac{ems} acts on a much larger time scale, a steady-state model for droop control is sufficient \cite{hans2021}. 
This model describes steady-state power sharing through the frequency-deviation variable $\rhok:=\omega^{\mathrm{nom}}-\omega(k)\in\RR$, where $\omega^{\mathrm{nom}}$ and $\omega(k)$ denote the nominal and steady-state frequencies, respectively, and the $\chi_i$ inverse droop gains.
In addition, physical power and energy limits are enforced through saturation in units. 
Let us now describe the saturation-based droop control model for different unit types based on \cite{HSH2021} using the notations $u_\nu\ok \in \mathbb{R}^{n_{\nu}}$ and $\chi_\nu \in \mathbb{R}^{n_{\nu}}$ for the power setpoint and inverse droop gain, of unit type $\nu\in\{\tth,\tst,\trs\}$ corresponding to conventional, storage, and renewable units. 

Conventional units, which usually serve as backup sources, may be switched on or off. Let us denote this switching behavior with the vector $\delta_\tth\ok \in \{0, 1\}^{\nt}$, 
where $\deltatik = 0$, indicates that conventional unit $i\in\{1,\dots, \nt\}$ is off, and $\deltatik = 1,$ indicates that it is on.
Thus, the saturated power output of the conventional unit is given by
\begin{equation}\label{eq:conventional:saturation:power}
\ptk = \text{diag}(\deltatk)\sat(\ptmin, \utk + \chi_{\tth}\,\rhok, \ptmax).
\end{equation}
Here, $\ptmin,\ptmax\in\Rpz^{\nt}$ are the conventional unit power limits and hence $\ptmin \leq p_t(k) \leq \ptmax$ is enforced through saturation for enabled units.

Storage units play a crucial role in dealing with the intermittency of renewable generation. Often, they are operated as grid-forming units. 
Let $\Ts \in \Rp$ be the sampling time, and $x(k) \in \Rpz^{\ns}$ be the stored energy level. Then, we have the storage dynamics,
\begin{subequations}
\label{eq:storage:model}
\begin{align}
x\ok &= x\okm - \Ts\,\psk,\\[-0.3em]
\intertext{and the storage energy limits}
\xsmin &\leq x\ok \leq \xsmax, \label{eq:storage:limits}
\end{align}
with $\xsmin \in \Rpz^{\ns}$ and $\xsmax \in \Rp^{\ns}$.

Storage operation is constrained by power and energy limits. To account for the energy limits, we define dynamically adjusted power limits based on the current energy and sampling time \cite{HSH2021}, \ie,
\label{eq:pstok}\begin{equation}\label{eq:psiminmaxk}\begin{split}
 \psmink &= \max\left(\psmin, (x\okm-x^\tmax)/\Ts\right),\\
 \psmaxk &= \min\left(\psmax, (x\okm-x^\tmin)/\Ts\right),
\end{split}\end{equation}
where $\psmin,\psmax\in\RR^{\ns}$ are the storage power limits.

Taking \eqref{eq:psiminmaxk} into account, the storage power is  
\begin{equation}\label{eq:pstokMain}
  \psk = \sat(\psmink, \usk+\chi_{\tst}\,\rhok, \psmaxk).
\end{equation}
\end{subequations}

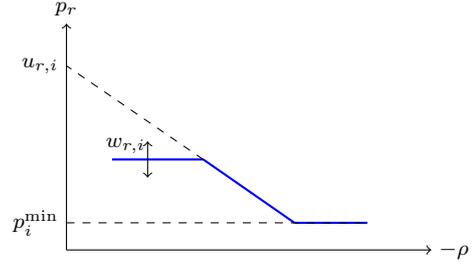
\begin{figure}
  \centering

\begin{tikzpicture}[scale=1.2]
   font=\footnotesize,
  \draw[->] (0,0) -- (4,0) node[right] {$-\rho = \omega - \omega^{\mathrm{nom}}$};
  \draw[->] (2,-0.3) -- (2,1.8) node[above] {$p_{\trs,i}$};

  \draw[very thick,draw=blue] (0.5,1) -- (1.5,1);          
  \draw[very thick,draw=blue] (1.5,1) -- (2.5,0);      
  \draw[very thick,draw=blue] (2.5,0) -- (3.3,0);
  \draw (2.15,0.4) -- (2.4,0.4);
  \draw (2.4,0.4) -- (2.4, 0.15);
  \node at (2.65,0.4) {$\chi_{\trs,i}$};
  \node[left] at (0,0) {$0$};
  \node at (2,-0.45) {$0$};
\node[left] at (2,0.5) {$u_{\trs,i}$};

  \draw[<->] (0.89,0.8) -- (0.89,1.2);
  \node at (0.65,1.15) {$\tiny w_{\trs,i}$};
\end{tikzpicture}
  \caption{Saturation-aware droop controlled RES unit $i$ with power setpoint
$u_{\trs,i}$ and available infeed $w_{\trs,i}$.}
\label{fig:resSat}
\end{figure}

We also assume that the \ac{res} have droop control. Still, their output power is limited by the available power under weather conditions $\wrk\in\Rpz^{\nr}$ and additionally at zero, \ie,
\begin{equation}
\label{eq:prsok}
 \prk = \sat(\mathbf{0}, \urk + \chi_\trs\,\rhok, \wrk).
\end{equation}
This is visualized in Figure~\ref{fig:resSat}. Note that $\wrk \leq \prmaxk$.


\section{Robust minimum-regret control problem}
\label{sec:problemFormulation}
We now formulate a multi-scenario \ac{mpc} problem that explicitly considers a range of possible scenarios from worst-case to best-case, and integrates  droop control and as well as physical limitations at all units through saturation. 
Our goal is to remain feasible for a range of admissible renewable and load realizations, and at the same time find feasible control actions whose cost is close to scenario-wise minima for all individual scenarios.

\subsection{Uncertainty model}
    \label{sec:uncertainModel}
Uncertainty is modeled in renewable generation and loads by lower and upper bounds. 
At each future time step $k$, a forecaster is assumed to provide these bounds such that

\begin{equation}
    \label{eq:uncertaintymodel}
    \underbrace{
    \begin{bmatrix}
    w_\trs^{\min}(k)\\
    w_\tld^{\min}(k)        
    \end{bmatrix}}_{w^{\min}(k)}
    \leq 
    \underbrace{
    \begin{bmatrix}
    w_\trs(k)\\
    w_\tld(k)        
    \end{bmatrix}}_{w(k)}
    \leq 
    \underbrace{
    \begin{bmatrix}
    w_\trs^{\max}(k)\\
    w_\tld^{\max}(k)        
    \end{bmatrix}}_{w^{\max}(k)},
\end{equation}
where $w_\trs^{\min}(k),w_\trs^{\max}(k) \in \mathbb{R}^{\nr}$ and $w_\tld^{\min}(k), w_\tld^{\max}(k) \in \mathbb{R}^{\nd}$ are the aforementioned upper and lower bounds of renewable generation and loads, respectively.

\subsection{Operating cost}
\label{sec:cost}
We consider an economically motivated operating cost associated with conventional and storage units, assuming zero operating cost for \ac{res}. 
Conventional units' cost includes fuel cost with weight $\Cth \in \Rpz^{\nt}$, fixed generation cost with weight $\CthOn\in \Rpz^{\nt}$ and switching cost with weight $\CthSw \in\Rpz^{\nt}$. 
Storage system cost has weight $\Cst \in \Rpz^{\ns}.$ Let us denote the power outputs and setpoints by 
  $p\ok = [p_\tth\ok\T ~ p_\tst\ok\T ~ p_\trs\ok\T]\T$ and, $u\ok = [u_\tth\ok\T ~ u_\tst\ok\T ~ u_\trs\ok\T ]\T,$ respectively. The total operation cost is 
    \begin{multline}\label{eq:stage:cost}
        \ell(p\ok, \deltatk, \deltath\okm) = \Cth\T\,\ptk + \CthOn\T\,\deltatk +\\
        \CthSw\T \left|\deltatk-\deltatkm\right| + \Cst\T\,\psk.
    \end{multline}
Note that $C_{\text{s},i}\psik$ is negative when charging, \ie, when $\psik < 0.$ 
This promotes storing available \ac{res} power that cannot be consumed immediately by the load.

\subsection{Multi-scenario minimum-regret operation}
Let $\Npr \in \N$ be the prediction horizon of the \ac{mpc}. At time instant $k,$ the predicted power at future time step $\tau \in\Npr$ is denoted by $p(k+\tau|k).$
From now, we follow the convention that if a vector $z\in \mathbb{R}^{n_z},$ then $\mathbf{z} = [z\okk~\cdots~z\okNp] \in \mathbb{R}^{n_z \times \Npr }.$
We define the operating cost as:
\begin{equation}\label{eq:operationalCost}
    J(\pVec, \deltathVec, \deltatk) := \sum_{\tau = 1}^{\Npr} \ell(p(k + \tau), \deltath(k+\tau), \deltath(k + \tau- 1)).
\end{equation}
This can be alternatively written as $J(\uVec, \deltathVec, \deltatk,\wVec)$ to highlight the dependence on the control input $\uVec$ and the uncertain input $\wVec.$
In the context of \ac{mpc}, the question arises as to whether there is a control trajectory that is not only feasible in every scenario (\ie, satisfies the model equations and constraints) but is also optimal in every scenario.

With the notation defined above, this can be expressed mathematically as the following robust \ac{mpc} problem.

\begin{problem}[Multi-scenario minimum-regret robust \ac{mpc}]%
\label{prob:multiSce}%
\begin{subequations}%
    \begin{equation}
        \label{eq:multimmmpc}
        \min\limits_{\uVec,\deltathVec} \max\limits_{\wVec} \big( J(\uVec, \deltathVec, \deltatk,\wVec) - \min\limits_{\uVec',\deltathVec'} J(\uVec', \deltathVec', \deltatk,\wVec) \big) 
    \end{equation}
subject to \eqref{eq:power:balance}--\eqref{eq:prsok} for $\uVec,\deltathVec$ and for $\uVec',\deltathVec'$. 
\begin{equation} \label{eq:minimax:control}       
    \deltathVec \in \{0, 1\}^{n_t\times \Npr}; \qquad
        \uVec^{\min} \leq \uVec \leq \uVec^{\max},  
    \end{equation} 
    where $\wVec\in[\wminVec,\wmaxVec]$ and $k$ is substituted by $k+\tau, \forall \tau \in \{1,\cdots,N_p\},$ with initial conditions 
    $x(k) = x^0$ and $ \deltath(k) = \delta_{\tth, 0}.$ Here $\uVec^{\min} \leq \uVec^{\max}$ are the control setpoint limits.
\end{subequations}
\end{problem}

Problem \ref{prob:multiSce} aims to minimize the cost across all possible values of the uncertainty variable $\wVec$. 
First $\min_{\uVec',\deltathVec'} J(\uVec', \deltathVec', \deltatk,\wVec),$ subject to the constraints of Problem~\ref{prob:multiSce}  gives the best possible cost for this specific scenario $\wVec.$ 
For each scenario, this gives us control trajectories $(\uVec'^\star,\deltathVec'^\star)$ which specifically minimize the cost for that scenario. We refer to this cost as scenario-specific optimal cost. 
In the second step, we search for one common control trajectory $(\uVec,\deltathVec)$ which is an optimal solution to Problem~\ref{prob:multiSce}. 
Let us consider a candidate $(\uVec,\deltathVec),$ and for each scenario, study the difference between the resulting cost $J(\uVec, \deltathVec, \deltatk,\wVec)$ and the previously calculated scenario-specific optimal cost. 
This difference is referred to as the regret, which tells us how much worse $(\uVec,\deltathVec)$ performs compared to $(\uVec'^\star,\deltathVec'^\star)$ for the same scenario $\wVec.$ 
In Problem~\ref{prob:multiSce}, we want to minimize the maximum of the difference between these two costs, \ie the maximum regret. 
In particular, $\max_{\wVec} $ chooses the disturbance that makes the regret worse and $\min_{\uVec,\deltathVec}$ aims to find the control trajectory $(\uVec^\star,\deltathVec^\star)$ that minimizes the worst-case regret. 
Clearly, if there exists a control $(\uVec^\star,\deltathVec^\star)$ which can optimize all scenarios simultaneously, then the optimal cost in \eqref{eq:multimmmpc} is zero.
Problem~\ref{prob:multiSce} can be challenging as it combines robust feasibility, binary variables, saturation, and an inner scenario-wise optimization. 
Therefore, in the next section, we develop an alternative tractable operation control approach based on a rule\hh based determination of constant power setpoints $\uVec$ and a unit\hh commitment problem\hh based determination of $\deltathVec$.


\section{Solution approach}
\label{sec:ruleBasedSolution}
In this section, we describe a tractable solution approach. 
First, we propose a simple rule-based controller under the assumption that all conventional units are always switched on. 
The resulting constant setpoints are subsequently used in a robust \ac{mpc}-based \ac{ems} which removes this assumption using a unit-commitment approach to decide which conventional unit is enabled or disabled. 

\subsection{Constant setpoint design with saturating droop control}
\label{sec:ruleBasedControl}

In what follows, we design a simple controller with constant setpoints. 
The design uses the saturation-based feasibility property discussed in \citep{HSH2021} which states that if power and energy saturation is enabled, the constraints \eqref{eq:power:balance}--\eqref{eq:prsok} are feasible if at any time, the power saturation limits allow that a power balance between generation and demand can be achieved for a given disturbance.
Furthermore, \cite{HSH2021} shows that if constraints \eqref{eq:power:balance}--\eqref{eq:prsok} are satisfied for $\wminVec$ and $\wmaxVec$, then they also hold for all $\wVec\in[\wminVec,\wmaxVec]$.
Let us first assume that all thermal generators are switched on $(\deltath=\mathbf{1})$, and that the load demand can be fully covered by the conventional units, \ie,
\begin{equation}
    \label{eq:conv:demand}
\mathbf{1}^T\ptmin \leq - \mathbf{1}^T w_d(k) \leq \mathbf{1}^T\ptmax.
\end{equation}

Under these assumptions, any given power setpoint trajectory $\uVec(k)\in [\uVec^{\min},\uVec^{\max}]$ is feasible, \ie, $\uVec^{\min}$ and $\uVec^{\max}$ can be chosen arbitrarily with $\uVec^{\min}<\uVec^{\max}$  because the feasibility depends on the saturation limits rather than on the setpoint values, and the always\hh on conventional units alone can provide sufficient power to meet the load demand.

In what follows, we assume that
\begin{enumerate}
    \item $\Cst\in[0,\Cth]$;
    \item the storage units can always be charged and discharged synchronously in accordance with their energy limits.
    In particular, let $\overline{\chi}_\tst$, $\overline{p}_\tst^\tmin$, $\overline{p}_\tst^\tmax$, $\overline{x}^\tmin$, $\overline{x}^\tmax$, $\overline{x}^0$ be the parameters and initial stored energy of a prototype storage unit. 
    The parameters and initial stored energy of each actual storage unit $i\in\{1,\ldots,n_\tst\}$  are obtained by scaling the parameters and initial stored energy of the prototype unit with a factor $\sigma_i$, \ie,
    \begin{multline*}
    	\chi_{\tst,i} = \sigma_i \overline{\chi}_\tst,\quad
    	p_{\tst,i}^\tmin = \sigma_i \overline{p}_\tst^\tmin,\quad
    	p_{\tst,i}^\tmax = \sigma_i \overline{p}_\tst^\tmax,\\
    	x_{\tst,i}^\tmax-x_{\tst,i}^\tmin = \sigma_i (\overline{x}^\tmax-\overline{x}^\tmin),\\
    	x_i^0-x_i^\tmin = \sigma_i (\overline{x}^0-\overline{x}^\tmin).
    \end{multline*}
    Then, the desired synchronized charging and discharging can be achieved when storage power setpoints are zero, \ie, $u_\tst = \mathbf{0}$.
\end{enumerate}
In addition, for simplicity, in what follows, we assume that all thermal units have the same cost factors $\Cth, \CthOn$ and all storage units have the same storage cost factor $\Cst$.

We now determine a set of \emph{constant} power setpoints $\uVec\in [\uVec^{\min},\uVec^{\max}]$ with the aim of minimizing the cost function \eqref{eq:operationalCost}.
This can be achieved by a rule\hh based controller, which is realized by appropriately selecting the limits and setpoints within the saturating droop\hh control laws.
In the end, we want to achieve the following behavior:
\begin{figure}[t]
        \centering

\begin{tikzpicture}
    font=\footnotesize,
\begin{axis}[
 width=0.5\textwidth, height=5cm,
  xmin=-4, xmax=4,
  ymin=-1.5, ymax=4,
  xtick={-4,-3,-2,-1,0,1,2,3,4},
  ytick={-1,0,1,2,3,4},
  xlabel={$-\rho = \omega - \omega^{\mathrm{nom}}$},
  ylabel={power [pu]},
  grid=both,
  major grid style={dotted, draw=black},
  legend cell align={left},
  legend style={at={(0.5, 1)}, anchor=south, draw=none, fill=none, column sep=0pt, },
  legend columns=3,
]

\addplot[very thick, green!70!black, mark=*, mark size=1.3pt]
  coordinates {(-4,1.0) (-3,1.0) (-2,1.0) (-1,0.2) (0,0.2) (1,0.2) (2,0.2) (3,0.2) (4,0.2)};
\addlegendentry{Conv. unit}

\addplot[very thick, red]
  coordinates {(-4,1.0) (-3,1.0) (-2,1.0) (-1,1.0) (0,0) (1,-1.0) (2,-1.0) (3,-1.0) (4,-1.0)};
\addlegendentry{Storage}

\addplot[very thick, orange, mark=triangle*, mark size=1.3pt]
  coordinates {(-4,1.2) (-3,1.2) (-2,1.2) (-1,1.2) (0,1.2) (1,1.2) (2,1.2) (3,0) (4,0)};
\addlegendentry{Wind turbine}

\addplot[very thick, blue]
  coordinates {(-4,0.55) (-3,0.55) (-2,0.55) (-1,0.55) (0,0.55) (1,0.55) (2.55,0.55) (3,0) (4,0)};
\addlegendentry{PV plant}

\addplot[ultra thick,densely dashdotted, black]
  coordinates {
    (-4,3.75) (-3,3.75) (-2,3.75)
    (-1,2.85) (0,1.90) (1,0.95) (2,0.95)
    (2.3,0.60) (2.6,0.20) (3,-0.80) (4,-0.80)
  };
  \addlegendentry{Sum of all units} 


\end{axis} 
\end{tikzpicture}
        \caption{Illustration of the rule-based configuration of the saturating droop functions.}
        \label{fig:rule-based-control}
\end{figure}
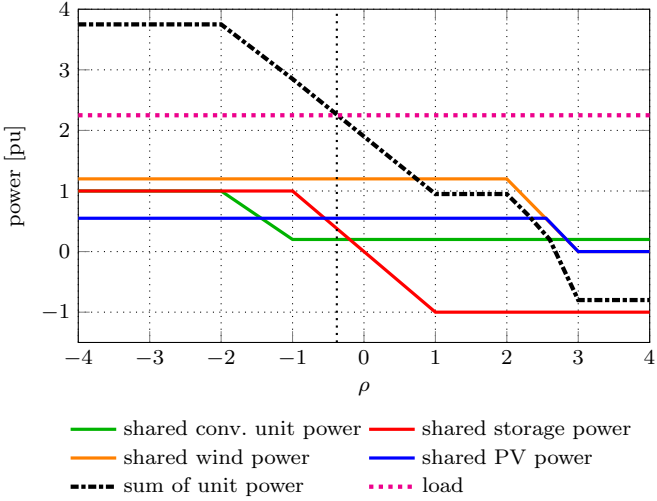
\begin{enumerate}
    \item use renewables (if available) and charge the battery (if not full),
    \item use renewables (if available) and discharge the battery (if not empty),
    \item use renewables (if available) and discharge the battery (if not empty) and use thermal generators.
\end{enumerate}
The saturating frequency\hh droop\hh curves for units are illustrated in Figure~\ref{fig:rule-based-control}.
For large $-\rho > 0$, all units are in lower saturation. For low $-\rho < 0$, units are in upper saturation.
For the above\hh mentioned prioritization to be realized by saturating droop curves, it is essential that the linear droop regions where $dp_{\nu}/d\rho \neq 0,$ of the unit types $\nu \in \{\text{t,s,r}\}$ are non-overlapping. 
This must be satisfied for all admissible saturation limits, and possible available renewable infeed since the upper saturation limit of a renewable unit depends on its available power, while the saturation limits of a storage unit depend on its \ac{soc}.
    
\begin{prop}\label{prop:ruleBasedControl}
Provided that $\umin$ is sufficiently low and $\umax$ is sufficiently high, possibly even outside of $[\pmin,\pmax]$, and all units are enabled, appropriate power setpoints $u^*$ can be calculated according to
    \begin{equation}
\label{eq:ruleBasedControl}
u^*_\tth=p^\tmin_\tth-\rho^\tmax_\tst\chi_\tth,\quad
u^*_\trs=p^\tmax_\trs-\rho^\tmin_\tst\chi_\trs,\quad
u^*_\tst=\mathbf{0},
\end{equation}
where $\rho^\tmin_\tst:=\min_{i \in [1,\ns] \subset N}(p^\tmin_{\tst,i}/\chi_{\tst,i})$ and
$\rho^\tmax_\tst:=\max_{i \in [1,\ns] \subset N}(p^\tmax_{\tst,i}/\chi_{\tst,i})$.
\end{prop}
\begin{proof} Let
\[
e_\nu(j\to N_p)
:=
T_s\sum_{l=j}^{N_p}\sum_{i=1}^{n_\nu}p_{\nu,i}(k+l),
\qquad \nu\in\{\text{t,r}\},
\]
denote the total infeed energy of unit type \(\nu\) from \(k+j\) to
\(k+\Npr\).
Since \(\deltath=\mathbf{1}\) is fixed, the on-cost due to $\CthOn$ is constant over the horizon and the switching cost $\CthSw$ is zero.
Furthermore, since all conventional and storage units have common cost factors $\Cth$ and $\Cst$, respectively, the cost reads
\begin{align*}
J
&=
\frac{\Cth}{\Ts}e_\tth(1\to \Npr)
+
\Cst\sum_{\tau=1}^{\Npr}
\mathbf{1}^T p_\tst(k+\tau)
+
\Npr\CthOn \nt.
\end{align*}
By summing the storage dynamics over the prediction horizon, we obtain
$\Ts\sum_{\tau=1}^{\Npr}\mathbf{1}^T p_\tst(k+\tau) = \mathbf{1}^T x(k)-\mathbf{1}^T x(k+\Npr).$ Therefore,
\begin{multline*}
J
=
\frac{\Cth}{\Ts}e_\tth(1\to \Npr)
+ 
\frac{\Cst}{\Ts}
\left(
\mathbf{1}^T x(k)-\mathbf{1}^T x(k+\Npr)
\right) \\
+
\Npr\CthOn \nt.
\end{multline*}
After removing the constant terms, let us shift the cost \eqref{eq:operationalCost} by a fixed part depending on initial stored energy and fixed on cost, \ie,
\begin{equation}
\begin{aligned}
\tilde J
&:=
J
-
\frac{\Cst}{\Ts}\mathbf{1}^T x(k)
-
\Npr\CthOn \nt \\
&=
\frac{\Cth}{\Ts}
\big(
e_\tth(\to \Npr)
-
\frac{\Cst}{\Cth}\mathbf{1}^T x(k+\Npr)
\big).
\end{aligned}
\label{eq:stage:costen}
\end{equation}

Consider the last time segment ($\Npr-1\rightarrow \Npr$), and the associated cost $\tilde{J}(\Npr-1\rightarrow \Npr)=\tfrac{\Cth}{\Ts}(e_\tth (\Npr-1\rightarrow \Npr) - \frac{\Cst}{\Cth}\mathbf{1}^T x\okNp)$.
The setpoints \eqref{eq:ruleBasedControl} make the non-saturated regions of
the saturating droop curves non-overlapping. Indeed, we have that
\[
p_\tth
=
\sat\!\left(
\ptmin,\,
\ptmin+\chi_\tth(\rho-\rho_\tst^\tmax),\,
\ptmax
\right).
\]
Hence, for \(\rho\leq\rho_\tst^\tmax\), the conventional units remain in lower saturation, i.e., \(p_\tth=\ptmin\). Similarly,
\[
p_\trs
=
\sat\!\left(
\mathbf{0},\,
\prmax+\chi_\trs(\rho-\rho_\tst^\tmin),\,
w_\trs
\right).
\]
So, the renewable units are at $w_\trs$ for $\rho\geq\rho_s^{\min}$. Together with \(u_s^*=0\), this results in the desired prioritization: renewable generation is used first, storage power is used next, and thermal generation only when required. 
Thus, for any fixed $w(k+N_p-1)$ and $x(k+N_p-1)$, \eqref{eq:ruleBasedControl} leads to minimal costs for time step $N_p-1\to N_p$.

In order to continue the proof by induction we need to know the effect that the control at time segment $\Npr-2\rightarrow \Npr-1$ has on the cost of the last time segment $\Npr-1 \rightarrow \Npr$. 
This impact is incurred due to how the control of segment $\Npr-2 \rightarrow \Npr-1$ changes the stored energy $x(k+\Npr-1)$ at the beginning of the last time segment. 
Since the optimal setpoints for the final time segment are given by \eqref{eq:ruleBasedControl} independently of $x(k+\Npr-1)$, it remains to study how a variation of this state affects the cost achieved under these setpoints.

Consider two states componentwise satisfying
\[
x^{(2)}(k+N_p-1)\geq x^{(1)}(k+N_p-1).
\]
For the same disturbance and the same control
\eqref{eq:ruleBasedControl}, the monotonicity properties of the
saturation-aware droop control derived in \cite{HSH2021}
imply that 
\[
p_\nu^{(2)}(k+N_p-1) \leq p_\nu^{(1)}(k+N_p-1),
\]
for $\nu\in\{t,r\},$ and
\[
x^{(2)}(k+N_p) \geq x^{(1)}(k+N_p).
\]
Consequently, for $\nu\in\{t,r\},$ we have
\[
e_\nu^{(2)}(N_p-1\to N_p) \leq e_\nu^{(1)}(N_p-1\to N_p).
\]
In particular, the larger stored energy at the beginning of the final segment cannot entail an increased optimal final\hh segment cost, and therefore
\begin{equation}
\tilde J^{(2)*}(N_p-1\to N_p)
\leq
\tilde J^{(1)*}(N_p-1\to N_p).
\label{eq:cost_monotonicity}
\end{equation}
    
Next, we also want to derive an upper bound of the cost reduction. For compactness, let us write for $\nu\in\{\text{t,r}\}$
\[
e_\nu^{(i)} := e_\nu^{(i)}(N_p-1\to N_p), \quad i\in\{1,2\}.
\]
Since both cases are evaluated for the same disturbance realization, the load demand over the final time segment is identical. For each
\(i\in\{1,2\}\), the energy balance over the final time segment is
\[ e_\tth^{(i)} + e_\trs^{(i)} + \Ts \mathbf{1}^T p_\tst^{(i)}(k+\Npr) +\Ts\,\mathbf{1}^T w_d(k+\Npr) =0.\]

Using the storage dynamics,
\[
\Ts\,\mathbf{1}^T p_\tst^{(i)}(k+\Npr)=\mathbf{1}^T x^{(i)}(k+\Npr-1)-\mathbf{1}^T x^{(i)}(k+\Npr),
\]
the energy balance for each case becomes
\begin{align*}
    e_\tth^{(i)} +e_\trs^{(i)}+&\mathbf{1}^T x^{(i)}(k+\Npr-1)\\
    &-\mathbf{1}^T x^{(i)}(k+\Npr)+ +\Ts\,\mathbf{1}^T w_d(k+\Npr) =0.
\end{align*}

Subtracting the balance equation for case~$2$ from case~$1$, yields
\begin{multline}
e_\tth^{(1)}-e_\tth^{(2)}
+
e_\trs^{(1)}-e_\trs^{(2)}\\
+
\mathbf{1}^T x^{(1)}(k+N_p-1)
-
\mathbf{1}^T x^{(2)}(k+N_p-1) \\
-
\left(
\mathbf{1}^T x^{(1)}(k+N_p)
-
\mathbf{1}^T x^{(2)}(k+N_p)
\right)
=0 .
\label{eq:energy_balance_two_cases}
\end{multline}
Hence, from \eqref{eq:stage:costen} and
\eqref{eq:energy_balance_two_cases},
\begin{multline}
\textstyle\frac{T_s}{C_t}
\big(
\tilde J^{(1)*}(N_p-1\to N_p)
-
\tilde J^{(2)*}(N_p-1\to N_p)
\big) = \\
\mathbf{1}^T x^{(2)}(k+N_p-1)
-
\mathbf{1}^T x^{(1)}(k+N_p-1) \\
-
\big(
1-\textstyle\frac{C_s}{C_t}
\big)
\big(
\mathbf{1}^T x^{(2)}(k+N_p)
-
\mathbf{1}^T x^{(1)}(k+N_p)
\big) \\
-\big(e_\trs^{(1)}-e_\trs^{(2)}\big).
\label{eq:cost_difference_bound_derivation}
\end{multline}
Since
\begin{align*}
    e_\trs^{(1)}-e_\trs^{(2)}&\geq0,\\
    \mathbf{1}^T x^{(2)}(k+N_p)-\mathbf{1}^T x^{(1)}(k+N_p)&\geq0,
\end{align*}
and \(C_s\leq C_t\), it follows from
\eqref{eq:cost_monotonicity} and
\eqref{eq:cost_difference_bound_derivation} that
\begin{equation}
\begin{aligned}
0
&\leq
\tilde J^{(1)*}(N_p-1\to N_p)
-
\tilde J^{(2)*}(N_p-1\to N_p) \\
&\leq
\textstyle\frac{C_t}{T_s}
\big(
\mathbf{1}^T x^{(2)}(k+N_p-1)
-
\mathbf{1}^T x^{(1)}(k+N_p-1)
\big).
\end{aligned}
\label{eq:storage_value_bound}
\end{equation}
Thus, if the total stored energy at the beginning of the final segment is
increased by
\[
\Delta x
:=
\mathbf{1}^T x^{(2)}(k+N_p-1)
-
\mathbf{1}^T x^{(1)}(k+N_p-1),
\]
then the corresponding reduction of the optimal remaining cost is at most
\((C_t/T_s)\Delta x\).

Now consider the preceding segment \(N_p-2\to N_p-1\). Any control different from \eqref{eq:ruleBasedControl} can affect the final-segment cost only through the stored energy at $k+\Npr-1$. 
If such a control increases the stored energy by $\Delta x,$ then by \eqref{eq:storage_value_bound}, the decrease in the final-segment cost is at most $(\Cth/\Ts)\Delta x$.
Due to prioritization in \eqref{eq:ruleBasedControl}, we already utilize available renewable infeed as much as possible.
Therefore, by power balance, it would require an increased generation of $\Delta x$ in conventional infeed over the preceding segment, with a cost increase of at least $(\Cth/\Ts)\Delta x$. 
Therefore, a deviation that increases the stored energy at $k+\Npr-1$ cannot reduce the total cost.
Thus, \eqref{eq:ruleBasedControl} is optimal for \(N_p-2\to N_p\).
Repeating the argument backward proves the claim over the entire prediction horizon by induction.
\end{proof}

\subsection{Extension by robust EMS with unit commitment}
\label{sec:unitCommitmentEms}
For calculating optimal constant power setpoints $u^*$ in the previous section, it was assumed that the thermal generators are always on.
The \ac{ems} problem defined next is not a direct solution to Problem~\ref{prob:multiSce}, but rather a tractable restricted formulation in which the continuous setpoints are fixed to the saturation-based values from Proposition~\ref{prop:ruleBasedControl} and the switch state of the generators are optimized using unit\hh commitment. 

The aim is to enable thermal generators only as long as necessary to robustly satisfy the constraints and minimize the worst\hh case cost $\max_{\wVec} J(u^*, \deltathVec, \deltatk,\wVec)$.
The associated unit commitment problem can be stated as follows.
\begin{problem}[Robust unit-commitment EMS]
        \label{prob:ruleBasedmpc}
    \begin{equation}
        \label{eq:ruleBasedmpc}
        \min_{\deltathVec}\ \max_{\wVec}\ J(u^*, \deltathVec, \deltatk, \wVec),
    \end{equation}
    subject to \eqref{eq:power:balance}--\eqref{eq:prsok} and constant setpoints \eqref{eq:ruleBasedControl} with initial conditions $x(k) = x^0,\quad \deltath(k) = \delta_{\tth, 0}$.
\end{problem}
Note that in contrast to Problem~\ref{prob:multiSce}, the continuous power setpoints are fixed in
Problem~\ref{prob:ruleBasedmpc}, and only the switching of the conventional units is optimized. 
Although Proposition~\ref{prop:ruleBasedControl} establishes the optimality of these setpoints when all conventional units remain switched on, this optimality is not directly claimed here for Problem~\ref{prob:ruleBasedmpc}, where conventional units may be enabled or disabled. 
Problem~\ref{prob:ruleBasedmpc} therefore constitutes a related tractable robust unit-commitment formulation rather than a direct solution of the minimum-regret Problem~\ref{prob:multiSce}.

\section{Case study}
\label{sec:caseStudy}
We evaluate the performance of Problem~\ref{prob:ruleBasedmpc}, on an islanded \ac{mg}. 
Therefore, we compare it with a prescient controller and the saturation-aware robust \ac{mpc} from \citep{HSH2021}.
\begin{figure}[t]
    \centering

\definecolor{mycolor1}{rgb}{0.06600,0.44300,0.74500}%
\definecolor{mycolor2}{rgb}{0.86600,0.32900,0.00000}%
\definecolor{mycolor3}{rgb}{0.92900,0.69400,0.12500}%
\definecolor{mycolor4}{rgb}{0.12941,0.12941,0.12941}%
\begin{tikzpicture}
  font = \footnotesize
\begin{axis}[%
width=0.4\textwidth, height=3cm,
at={(0cm,0cm)},
scale only axis,
xmin=0,
xmax=24,
xtick={0,4,8,12,16,20,24},
xlabel style={font=\color{mycolor4}},
xlabel={Time (in hours)},
ylabel ={Load power [pu]},
ymin=-1,
ymax=-0.3,
grid=both,
major grid style={dotted, draw=black},
legend columns=3,
legend style={
  at={(0.58,0.98)},
  anchor=north,
  legend cell align=left,
  align=left,
  draw=none,
  fill=none
}
]

\addplot [very thick, color=mycolor2]
table[row sep=crcr]{%
0	-0.741024125753995\\
1	-0.64682020562738\\
2	-0.403195676917048\\
3	-0.37810906551253\\
4	-0.549534340406625\\
5	-0.359210392011127\\
6	-0.564086161309377\\
7	-0.472879245285666\\
8	-0.624081708528562\\
9	-0.799064911950366\\
10	-0.580655940023691\\
11	-0.452024141918753\\
12	-0.559590588313557\\
13	-0.76254743900028\\
14	-0.497315326883765\\
15	-0.603843892353784\\
16	-0.603630618865215\\
17	-0.554122668959519\\
18	-0.510404594324234\\
19	-0.780316061480191\\
20	-0.521375153262009\\
21	-0.534068267577997\\
22	-0.633175481080749\\
23	-0.767521240917463\\
24	-0.692382164609842\\
};
\addlegendentry{$\alpha = 1$}


\addplot [very thick, color=mycolor3]
 table[row sep=crcr]{%
0	-0.816024125753995\\
1	-0.72182020562738\\
2	-0.478195676917048\\
3	-0.45310906551253\\
4	-0.624534340406625\\
5	-0.434210392011127\\
6	-0.639086161309377\\
7	-0.547879245285666\\
8	-0.699081708528562\\
9	-0.874064911950366\\
10	-0.655655940023691\\
11	-0.527024141918753\\
12	-0.634590588313558\\
13	-0.83754743900028\\
14	-0.572315326883765\\
15	-0.678843892353784\\
16	-0.678630618865215\\
17	-0.629122668959519\\
18	-0.585404594324234\\
19	-0.855316061480191\\
20	-0.596375153262009\\
21	-0.609068267577997\\
22	-0.708175481080749\\
23	-0.842521240917463\\
24	-0.767382164609842\\
};
\addlegendentry{$\alpha = 0.5$}

\addplot [very thick, color=mycolor1]
  table[row sep=crcr]{%
0	-0.891024125753995\\
1	-0.79682020562738\\
2	-0.553195676917048\\
3	-0.52810906551253\\
4	-0.699534340406625\\
5	-0.509210392011127\\
6	-0.714086161309377\\
7	-0.622879245285666\\
8	-0.774081708528562\\
9	-0.949064911950366\\
10	-0.730655940023691\\
11	-0.602024141918753\\
12	-0.709590588313558\\
13	-0.91254743900028\\
14	-0.647315326883765\\
15	-0.753843892353784\\
16	-0.753630618865216\\
17	-0.704122668959519\\
18	-0.660404594324234\\
19	-0.930316061480191\\
20	-0.671375153262009\\
21	-0.684068267577997\\
22	-0.783175481080749\\
23	-0.917521240917463\\
24	-0.842382164609842\\
};
\addlegendentry{$\alpha = 0$}
\end{axis}
\end{tikzpicture}%
    \caption{Load demand profiles $\wld^{\alpha}$ for selected scenarios on day~$1$.}
    \label{fig:loadScenario}
\end{figure}
\subsection{Simulation setup}
\label{simulationSetup}
We consider an islanded microgrid comprising a conventional generator, a storage unit, a PV plant, a wind turbine and a load.
The unit parameters, cost coefficients, and initial conditions are listed in Table~\ref{tab:caseStudy}.
\begin{table}[h]
\centering
\caption{Unit parameters, cost coefficients, and initial conditions for the case study.}
\label{tab:caseStudy}
\begin{tabular}{ll@{\hspace{1.2em}}ll}
\toprule
Parameter & Value & Parameter & Value \\
\midrule
$\ptmin$ & $0.2\,\unit{pu}$ & $\Cth$ & $1$ \\
$\psmin$ & $-1\,\unit{pu}$ & $\CthOn$ & $0.2$ \\
$\ptmax$ & $1\,\unit{pu}$ & $\CthSw$ & $0.3$ \\
$\psmax$ & $1\,\unit{pu}$ & $\Cst$ & $0.9$ \\
$x_\tmin$ & $0\,\unit{pu\,h}$ & $\delta_{\tth,0}$ & $0$ \\
$x_\tmax$ & $6\,\unit{pu\,h}$ & $x^0$ & $2\,\unit{pu\,h}$ \\
$u_i$ & $[-5,5]\,\unit{pu}$ & $\chi_i$ & $1$ \\
\bottomrule
\end{tabular}
\end{table}
We choose $\Cst<\Cth$ to penalize charging from conventional generation.

The sampling time is \(15\,\unit{min}\), \(N_p=32\), and the closed-loop simulation length is \(\Nsi=672\), corresponding to \(7\,\unit{days}\). 
Simulations use Matlab\textsuperscript{\textregistered} 2022b, YALMIP~\citep{Lof2004}, and Gurobi~11.0.0~\citep{gurobi}.

\subsection{Prescient controller}\label{sec:prescientMpc}
A prescient \ac{mpc} controller which is based on Problem~\ref{prob:multiSce} is used as a benchmark. It is assumed to know the realized available renewable power and the demand trajectory over the prediction horizon, denoted by $\mathbf w^{\text{actual}}$. 
Therefore, the uncertainty set in \eqref{eq:uncertaintymodel} is substituted by $\wminVec =\wmaxVec = \mathbf{w}^{\text{actual}},$ effectively reducing the multi-scenario \ac{mpc} in \eqref{eq:multimmmpc} to a single-scenario problem. 

\subsection{Closed-loop simulations}
We compare the performance of different controllers over different disturbance scenarios. 
Since \(w_d\leq 0\), the worst-case scenario $\wminVec$ corresponds to minimum available renewable power and maximum absolute value of power consumption, whereas the best-case scenario $\wmaxVec$ corresponds to maximum renewable available power and minimum absolute value power consumption. 
We create $11$ scenarios using linear interpolation between the worst-case and the best-case scenarios via
\begin{align}\label{eq:disturbanceScenario}
    w^{\alpha}(k) = w^{\min}(k) + \alpha \big(w^{\max}(k) - w^{\min}(k)\big), \quad \forall k
\end{align}
where $\alpha\in \{0, 0.1,\dots, 1\}$. Hence, $\alpha = 0,$ and $\alpha=1$ correspond to the worst-case and best-case scenarios, respectively.
As an example, Figure~\ref{fig:loadScenario} shows different load-demand scenarios $\wld^{\alpha}.$ The total cost over the entire simulation horizon of $\Nsi=672$ samples is computed from the stage cost \eqref{eq:stage:cost} as
\begin{equation}
 J^\text{closed\hh loop}:=\textstyle\sum_{k=1}^{\Nsi}\ell\left(p\ok,\deltath\ok,\deltath\okm\right),
\end{equation}
with the model parameters listed in Table~\ref{tab:caseStudy}.

\begin{figure}[t]

\definecolor{mycolor1}{rgb}{0.06600,0.44300,0.74500}%
\definecolor{mycolor2}{rgb}{0.86600,0.32900,0.00000}%
\definecolor{mycolor3}{rgb}{0.92900,0.69400,0.12500}%
\definecolor{mycolor4}{rgb}{0.12941,0.12941,0.12941}%
\begin{tikzpicture}
  font = \footnotesize
\begin{axis}[%
width=0.41\textwidth, height=4cm,
at={(0cm,0cm)},
scale only axis,
xmin=1,
xmax=11,
xtick={1,3,5,7,9,11},
xticklabels={$0$,$0.2$,$0.4$,$0.6$,$0.8$,$1$},
xlabel style={font=\color{mycolor4}},
xlabel={$\alpha$},
ymin=-50,
ymax=350,
ylabel style={font=\color{mycolor4}},
ylabel={Cost},
grid=both,
major grid style={dotted, draw=black},
legend cell align={left},
legend style={
  at={(0.4,0.97)},
  anchor=north west,
  draw=none,
  fill=white,
  fill opacity=0.05,
  text opacity=1,
  font=\scriptsize
},
legend columns=1
]

\addplot [very thick, color=black]
  table[row sep=crcr]{%
1	307.384414334523\\
2	262.305185509041\\
3	223.77343801185\\
4	186.913462663979\\
5	141.397402605771\\
6	90.6809627569409\\
7	38.2741833623943\\
8	0.893657605972165\\
9	-2.27114583386379\\
10	-5.84642295557897\\
11	-8.94931020852474\\
};
\addlegendentry{Prescient controller}

\addplot [very thick, color=mycolor1]
  table[row sep=crcr]{%
1	314.354234962196\\
2	276.581118433685\\
3	236.259468831966\\
4	194.212756137072\\
5	150.105695491651\\
6	98.5767723917597\\
7	42.3954936122989\\
8	3.82937653012716\\
9	-0.155339176297863\\
10	-3.68301648435637\\
11	-6.92101501332761\\
};
\addlegendentry{Robust unit commitment EMS}

\addplot [very thick, color=mycolor2]
  table[row sep=crcr]{%
1	299.484414331037\\
2	269.327442877423\\
3	242.351299721647\\
4	224.239640752009\\
5	201.56582059775\\
6	172.83694677959\\
7	152.981700816641\\
8	129.217511918502\\
9	100.247537602647\\
10	71.6572455865193\\
11	57.0369124824268\\
};
 \addlegendentry{Sat.\hh aware robust MPC \cite{HSH2021}}
\end{axis}
\end{tikzpicture}%
    \caption{Comparison of the closed-loop cost calculated over one week, for different controllers and disturbance scenarios.}
    \label{fig:costComparison}
\end{figure}

Figure~\ref{fig:costComparison} depicts the total closed-loop cost associated with the considered controllers, for each scenario defined in \eqref{eq:disturbanceScenario}. 
The proposed robust unit-commitment \ac{ems} in Problem~\ref{prob:ruleBasedmpc}, achieves a cost  close to that of the prescient controller for all scenarios. 
Moreover, it achieves a lower closed-loop cost than the saturation-aware minimax MPC \citep{HSH2021}, for non-worst-case scenarios. 
This improvement is consistent with the use of the priority-based constant power setpoints, whose optimality is established under the always-on assumption of Proposition~\ref{prop:ruleBasedControl}, together with the robust switching of the conventional generator in Problem~\ref{prob:ruleBasedmpc}. 
Under less adverse disturbance scenarios, these priority\hh based constant setpoints exploit the additional available renewable power and lower demand through saturation\hh based power sharing.

The observed slight performance gap could be attributed to the restriction to fixed priority\hh based constant setpoints in the presence of switching of conventional units. 
When fixed on cost and switching cost are relevant, it may be beneficial to operate the conventional more intensively over shorter on-periods, store the generated energy, and use it later. 
Such behavior cannot generally be represented by fixed setpoints of Proposition~\ref{prop:ruleBasedControl}, whose optimality is established under the always-on assumption. 
However, overall the approach performs well with a cost very close to that of the prescient \ac{mpc}.


\section{Conclusion}
\label{sec:conclusion}
In this paper, we presented a novel hierarchical control strategy for islanded microgrids.
On the primary control layer, droop control of conventional, storage and renewable units is considered, which also includes autonomous power and energy limiting. 
We have shown that appropriately chosen constant power setpoints can lead to an optimal control when assuming always-on conventional generator. 
This property is then used in a robust MPC-based energy management system that decides which conventional units are enabled or disabled to minimize operating cost under a wide range of scenarios of uncertain inputs. 
In a case study, we demonstrated that the proposed strategy achieves in almost all cases achieves lower costs than existing robust energy management schemes. 
Given the simplicity of our approach, we suggest that it can be especially well suited for small-scale microgrids. 
Future work will focus on larger microgrids and broader uncertainty sets.
\bibliographystyle{IEEEtran}
\bibliography{ieeeconf,rulebased_paper}
\end{document}